\documentclass{svjour3}
\smartqed  
\usepackage{graphicx}
\usepackage{latexsym,bm,amsmath,amssymb}
\usepackage{rotating}
\usepackage{multirow,booktabs,cases}
\usepackage[misc]{ifsym}
\usepackage[style=1]{mdframed}
\usepackage{algorithm,algorithmic}
\usepackage[colorlinks=true]{hyperref}
\hypersetup{urlcolor=blue, citecolor=blue,linkcolor=blue}
\usepackage{epstopdf}
\usepackage{bm}

\def\x{{\bf x}}
\def\y{{\bf y}}
\def\z{{\bf z}}

\def\a{{\bf a}}
\def\b{{\bf b}}

\def\0{{\bf 0}}
\def\A{{\cal A}}
\def\B{{\cal B}}
\def\C{{\cal C}}

\def\X{{\cal X}}

\def\Y{{\cal Y}}
\allowdisplaybreaks

\begin{document}
\graphicspath{{./PIC/}}

\title{Completely Positive Reformulations of Polynomial Optimization Problems with Linear Inequality Constraints}
\titlerunning{Completely Positive Reformulations of POP with LIC}     

\author{Haibin Chen \and Hong Yan \and Guanglu Zhou}


\institute{Authors are listed alphabetically. \\
H. Chen \at School of Management Science, Qufu Normal University, Rizhao, Shandong, China, 276800. \\
\email{chenhaibin508@163.com}
\and H. Yan \at Centre for Intelligent Multidimensional Data Analysis, Hong Kong Science Park, Shatin, Hong Kong. Department of Electrical Engineering, City University of Hong Kong, Hong Kong, China.\\
\email{h.yan@cityu.edu.hk}
\and G. Zhou \at Curtin Centre for Optimisation and Decision Science,
School of Electrical Engineering, Computing and Mathematical Sciences (EECMS),
Curtin University, Perth, WA 6845, Australia. \\
\email{g.zhou@curtin.edu.au}}

\date{Received: date / Accepted: date}

\maketitle

\begin{abstract}
Polynomial optimization encompasses a broad class of problems in which both the objective function and constraints are polynomial functions of the decision variables. In recent years, a substantial body of research has focused on reformulating polynomial optimization problems (POPs) as conic programs over the cone of completely positive tensors (CPTs).  In this article, we propose several new completely positive reformulations for a class of POPs with linear inequality constraints. Our approach begins by lifting these problems into a novel convex optimization framework, wherein the variables are represented as combinations of symmetric rank-one tensors. Based on this lifted formulation, we present a general characterization of POPs with linear inequality constraints that can be reformulated as conic programs over the CPT cone. Additionally, we construct the dual formulations of the resulting completely positive programs. Under mild assumptions, we prove that these dual problems are strictly feasible and strong duality holds.

\end{abstract}

\section{Introduction}\label{Int}
The polynomial optimization problem (POP) is an optimization model in which objective and constraints are
polynomial functions. It has attracted much research attention in the past decades. This is in part due to its own strong theoretical appeal
and an increased demand on applications in
numerical linear algebra \cite{CCL,Hillar2013,Qi2005}, material sciences \cite{HC18,SYC2008}, quantum physics \cite{DLMO07,ZNZ20},
and signal processing \cite{QT03,san2022}.
Generally, solving the polynomial optimization problem is a challenging task as it is NP-hard. For instance, even the simplest instances of polynomial optimization, such as maximizing a cubic polynomial over a sphere, are NP-hard \cite{Nest03,So11}.
For readers interested in polynomial optimization with
simple constraints, see \cite{Klerk08} for a survey on the computational complexity
of optimizing various classes of polynomial functions over a simplex, hypercube, or
sphere.

In mathematics, there are two main research directions for POP. The first direction aims to get an approximate solution of POP by algorithms designed directly for the given POP. The second direction pays attention to getting near optimal solutions or lower (or upper) bounds for the optimal value of the original POP by its relaxation problems. In the algorithm point of view, researchers are encouraged to develop some effective approximation algorithms for finding approximately optimal solutions, which include semi-definite programming (SDP) method \cite{HLZ10,JB06,JB09,Ling17,Ling09,LZ10,Nie13}, alternating direction method of multipliers \cite{Jiang}, block coordinate descending method \cite{CHZW,RHL13,Tseng01} and some variations of the algorithms above \cite{Chen2012,Zhou12,Nest97,shor87}.
On the other hand, a well established line of research is that the POP can be reformulated (or relaxed) as the completely positive tensor program (CPTP) or its dual problem i.e., the copositive tensor program (CTP).
The CPTP and CTP are linear programs over the cone of completely positive tensor cone and its dual cone i.e., the cone of copositive tensors. A pioneer contribution to the line of research can be found in \cite{pena14}, and more details about these two kinds of cones can be found in  \cite{HC18,CH17,zf14}.

In 2015, Pe$\tilde{\rm n}$a {\it et al.} in \cite{pena14} provided a general characterization of the class of polynomial optimization problems that can be
formulated as a conic program over the cone of completely positive tensors. As a consequence of
this characterization, it follows that recent related results for quadratic problems can be further
strengthened and generalized to high-order polynomial optimization problems. Zuluaga {\it et al.} in \cite{Kuang18}
studied CPTP relaxations and completely positive semidefinite tensor relaxations for general POPs under special constraints. By these relaxations, they
gave tighter bounds than the bounds obtained via a Lagrangian relaxation of the POPs. Later, Xia {\it et al.} in \cite{XIAZ2017} proposed several novel convex reformulation results for POPs which consider the set
of linear constraints explicitly, where the results can be viewed as a refinement of Theorem 5 in \cite{pena14}, and the generalization of Theorem 4 in \cite{Bai2016}.
Furthermore, the CPTP method can be applied to practical problems. An important example of this approach is the work of Chen {\it et al.} \cite{CHZW} in which they showed that
computing the coclique number of a uniform hypergraph can be reformulated as a CPTP.

In this paper,
we consider the equivalent CPTP reformulations of POPs, where instead of assuming the polynomial constraints to be nonnegative over $\mathbb{R}^n_+$ \cite{pena14}, we only assume the nonnegativity of polynomial constraints over the feasible set of linear inequality constraints.
As a result of this reformulation approach, it is desirable to establish a deeper connection between polynomial optimization and conic programming. This connection enables the development of novel solution schemes for general polynomial optimization problems by leveraging tools from conic programming, particularly the extensively studied cones of completely positive and copositive tensors.
The main distinction between the results presented in this paper and those in \cite{pena14} lies in the treatment of inequality constraints: while \cite{pena14} introduces slack variables to handle inequalities, our approach directly addresses both equality and inequality constraints without the need for such transformations.
To this end, we begin by introducing a general lifting framework that maps a POP into a convex program involving tensor variables. By analyzing the recession cone of the corresponding feasible set, we demonstrate that all feasible POPs can be equivalently reformulated as convex problems whose variables are combinations of symmetric rank-one tensors.
Building on this framework, we provide a general characterization of POPs with linear inequality constraints that can be reformulated as CPTPs. We also construct the dual formulations of the resulting CPTPs, which take the form of copositive tensor programs. We show that these dual problems are strictly feasible and that strong duality holds under mild conditions.
Consequently, existing numerical methods developed for both the primal and dual problems can be effectively applied to solve these reformulated models.

The remainder of this paper is organized as follows. In Section 2, we review some basic definitions related to tensors and polynomials. Section 3 introduces the lifting model for general POPs and establishes its equivalence to the original formulation. In Section 4, we present the equivalent completely positive reformulations for POPs with linear inequality constraints and demonstrate that strong duality holds between the resulting CPTPs and their duals. Final remarks and conclusions are provided in Section 5.

\section{Preliminaries}\label{Sec2}
In this section, we first review the necessary notation and relevant preliminaries. We then provide a detailed comparison between our results and the existing reformulations presented in \cite{pena14,XIAZ2017}.

\subsection{Basic notation}
For a positive integer $n$, denote $[n]=\{1,2,\cdots,n\}$.
Real scalars from $\mathbb{R}$ are denoted by unbolded lowercase letters such as $\lambda$, $a$ and denote vectors by bolded lowercase letters
such as $\x$, $\y$. Denote $\mathbb{R}^n$($\mathbb{R}^n_+$) the set including all $n$ dimensional vectors(nonnegative vectors), and let $\mathbb{R}^{m\times n}$ be the set including all $m\times n$ matrices.
Furthermore, matrices are denoted by uppercase letters such as $A, B$. Tensors are denoted by capital uppercase letters such as $\A$, $\B$. An order $d$ dimension $n$ tensor $\A$ can be defined by $\A=(a_{i_1i_2\cdots i_d})$,
where $i_1,\cdots,i_d\in[n]$ and $a_{i_1i_2\cdots i_d}$ is the $(i_1,\cdots,i_d)$-th entry of $\A$. Tensor $\A$ is called symmetric if
$a_{i_1i_2\cdots i_d}=a_{i_{\sigma(1)}i_{\sigma(2)}\cdots i_{\sigma(d)}}$ for any permutation $\sigma$ of $[d]$.
We denote $S^{d,n}$ as the set of all symmetric tensors of order $d$ and dimension $n$.

Let $\A=(a_{i_1i_2\cdots i_d})$, $\B=(b_{i_1i_2\cdots i_d})\in S^{d,n}$ be two symmetric tensors. The inner product of $\A$ and $\B$ is defined such as
$$
\langle\A,~ \B\rangle=\sum_{i_1,i_2,\cdots,i_d\in[n]}a_{i_1i_2\cdots i_d}b_{i_1i_2\cdots i_d}.
$$
Suppose that $D_1, D_2,\cdots, D_d\in\mathbb{R}^{p\times n}$. The mode production between tensor $\A$ and the matrices is denoted by
$\A_{\times 1}{D_1}_{\times 2}{D_2}\cdots_{\times d}D_d$ with entries such that
$$
(\A_{\times 1}{D_1}_{\times 2}D_2\cdots_{\times d}D_d)_{i_1i_2\cdots i_d}=\sum_{j_1,j_2,\cdots j_d\in[n]}a_{j_1j_2\cdots j_d}(D_1)_{i_1j_1}\cdots(D_d)_{i_dj_d},
$$
for all $i_1, i_2, \cdots, i_d\in[p]$. Note that $\A_{\times 1}{D_1}_{\times 2}{D_2}\cdots_{\times d}D_d$ may not be symmetric.
Moreover, for $\X\in\mathbb{S}^{d,p}$ it holds that
$$
\langle\X, ~\A_{\times 1}{D_1}_{\times 2}D_2\cdots_{\times d}D_d\rangle=\langle\X_{\times 1}{D^\top_1}_{\times 2}D^\top_2\cdots_{\times d}D^\top_d, ~\A\rangle.
$$

Let $\mathbb{R}_d[\x]$ denote the set of real coefficient polynomials of degree at most $d$.
It is well known that any given high order homogeneous polynomial corresponds uniquely to a symmetric tensor.
To move on, we introduce the following definition, which maps a vector to a symmetric rank-1 tensor.
\begin{definition}\label{def32}{\rm\cite{pena14}} For any $\x\in\mathbb{R}^n$, denote the mapping $M_d: \mathbb{R}^n\rightarrow S^{d,n}$ such that
$$
M_d(\x)=\underbrace{\x\circ\x\circ\cdots\circ\x}_d,
$$
where $``\circ"$ denotes the out product and $\underbrace{\x\circ\x\circ\cdots\circ\x}_d$ denotes a symmetric rank-1 tensor.
\end{definition}

For a homogeneous polynomial $f(\x)\in \mathbb{R}_d[\x], \x\in\mathbb{R}^n$, its corresponding tensor is denoted by $\A_f=(a_{i_1i_2\cdots i_d})\in S^{d, n}$, then we have the following representation:
$$
f(\x)=\A_f\x^d=\langle\A_f, M_d(\x)\rangle=\sum_{i_1,\cdots,i_d\in[n]}a_{i_1i_2\cdots i_d}x_{i_1}x_{i_2}\cdots x_{i_d}.
$$
In Section 3, we establish a mapping from polynomial systems to the set of symmetric tensors and demonstrate that non-homogeneous polynomials retain a one-to-one correspondence with symmetric tensors.

For any subset $F\subseteq\mathbb{R}^n$, the recession cone of $F$ is defined as follows:
$$
REC_F=\{\y\in\mathbb{R}^n~|~\forall~ \x\in F,~\forall~\lambda\geq 0,~\lambda\in\mathbb{R},~ \x+\lambda\y\in F \}.
$$
Moreover, the convex sets below are useful in the sequel:
\begin{equation}\label{e21}
C^h_F={\rm conv}\{M_d(\x)~|~\x\in F\},~~R^h_F={\rm conv}\{M_d(\x)~|~\x\in REC_F\}.
\end{equation}
By the notion of $REC_F$, it is clear that $R_F$ is a convex cone. The order of the tensors in $C^h_F$ and $R^h_F$ are clear through the text.
In general, computing the recession cone for high-degree polynomial constraints is challenging, though tractable in certain special cases \cite{KU2023}. The following example illustrates recession cones corresponding to simple feasible sets defined by linear and quadratic constraints, respectively.
\begin{example}\label{exam21}Suppose that $A\in\mathbb{R}^{m\times n}$ is a matrix and $\b\in\mathbb{R}^m$. If $F=\{\x\in\mathbb{R}^n~|~A\x+\b\geq\0\}$, then $REC_F=\{\y\in\mathbb{R}^n~|~A\y\geq\0\}$. Furthermore, assume $B\in\mathbb{R}^{n\times n}$ is a symmetric matrix. If $F=\{\x\in\mathbb{R}^n~|~\x^\top B\x=0 \}$, then
$REC_F=\{\y\in\mathbb{R}^n~|~\y^\top B\y=0,~\x^\top B\y=0,~\forall~\x\in F\}$.
\end{example}

To end this section, we recall definitions of the copositive tenor and the completely positive tensor. For $\A, \X\in S^{d,n}$, $\A$ is called a copositive tensor(or called copositive for simple) if $\A\x^d\geq 0,~{\rm for~ all}~\x\in\mathbb{R}^n_+$, and $\X$ is a completely positive tensor if it can be decomposed into the sum of several nonnegative symmetric rank one tensors i.e., $\X=\sum_{i=1}^{k}M_d(\x_i)\in\mathbb{R}^n_+$.
The copositive tensor cone and the completely positive tensor cone are denoted by $\mathcal{C}_{n,d}$ and $\mathcal{C}^*_{n,d}$ respectively. It is well known  that $\mathcal{C}_{n,d}$ and $\mathcal{C}^*_{n,d}$ are dual to each other.

\subsection{Completely positive tensor reformulations on POPs}
It is known that under appropriate conditions, a POP can be reformulated as a
CPTP, that is, a conic (linear) program over the cone of
completely positive tensors (see, e.g., \cite{pena14,XIAZ2017}). Next, we summarize the theorems that provide insights on the conditions that guarantee the equivalence of the POPs with their respective CPTP relaxations.

In the pioneering work \cite{pena14}, Pe$\tilde{\mbox n}$a et al. considered the following POPs:
\begin{equation}\label{e22-1}
\inf ~q(\x)~~~\mbox {s.t.}~ h_i(\x)=0, \x\geq\0, i=1,\ldots,m,
\end{equation}
for some given $n$-variate polynomials $q, h_i\in\mathbb R_d[\x]$. The following CPTP ia a relaxation of (\ref{e22-1}):
\begin{equation}\label{e22-2}
\begin{array}{ll}
\inf~&\langle C_d(q), \Y\rangle\\
\mbox{s.t}~&\langle C_d(h_i), \Y\rangle=0, i=1,2,\ldots,m,\\
&\langle C_d(1), \Y\rangle=1,\\
&\Y\in\mathcal{C}^*_{n+1,d},
\end{array}
\end{equation}
where $C_d(q)$ is the corresponding coefficient tensor of polynomial $q(\x)$.

Note that (\ref{e22-2}) is a natural convex lifting of the POP (\ref{e22-1}), and the relaxation is not always tight as illustrated in \cite{pena14}. In order to character conditions under which the relaxation (\ref{e22-2}) is tight, further verification are needed. We summarize the sufficient conditions in Theorem \ref{th-1}, which were first raised in \cite{pena14}. To move on, let $\tilde{h}(\x)$ denote the homogeneous component of $h$ of highest degree. For nonempty set $S\subseteq\mathbb{R}^n$, its horizon cone is defined as
$$
S^{\infty}:=\{\y\in\mathbb{R}^n : \exists~\x^{k}\in S, \lambda^k\geq 0, \lambda^k\in\mathbb{R}, k=1,2,\ldots, \mbox{such that}~ \lambda^k\rightarrow 0, \lambda^k\x^k\rightarrow\y\}
$$

\begin{theorem}\label{th-1}{\rm\cite{pena14}} Let $q, h_i\in\mathbb{R}_d[\x]$ for $i=1,\ldots,m$. Then,  problems (\ref{e22-1}) and (\ref{e22-2}) are equivalent in any of the following two cases.

\noindent(i) $\deg(h_i)=d, h_i(\x)\geq0$ for all $\x\in S_{i-1}$, and $\{\x\in S_{i-1}^{\infty} : \tilde{h}_i(\x)=0\}\subseteq S_i^{\infty}$, where $S_0=\mathbb{R}^n_+$ and $S_i=\{\x\in S_{i-1} : h_i(\x)=0\}, i=1,\ldots,m$.

\noindent(ii) $\deg(h_i)=d, h_i(\x)\geq0$ for all $\x\in\mathbb{R}^n_+$, and $\tilde{q(\x)}\geq 0$ for all $\x\in\{\x\in\mathbb{R}^n_+ : \tilde{h}_i(\x)=0, i=1,\ldots,m\}$.
\end{theorem}

Xia et al. \cite{XIAZ2017} later refined the sufficient condition $(ii)$ presented in Theorem \ref{th-1}, and considered the following POP:
\begin{equation}\label{e22-3}
\begin{array}{ll}
\inf~&q(\x) \\
\mbox{s.t.}~&l_i(\x)=0, i=1,\ldots,m_l\\
&h_j(\x)=0,j=1,\ldots,m_n\\
&\x\geq\0,
\end{array}
\end{equation}
where $l_i(\x)=\a_i^\top\x-b_i=0$ with $\a_i\in\mathbb{R}^n$ and $b_i\in\mathbb{R}$. 
Accordingly, the CPTP relaxation is formulated as follows:
\begin{equation}\label{e22-4}
\begin{array}{ll}
\inf~&\langle C_d(q), \Y\rangle \\
\mbox{s.t.}~&\langle C_d(h_j), \Y\rangle=0, j=1,\ldots,m_n\\
&\langle C_d(l_i^d), \Y\rangle=0, i=1,\ldots,m_l\\
&\langle C_d(1), \Y\rangle=1\\
&\Y\in\mathcal{C}^*_{n+1,d}.
\end{array}
\end{equation}

To prove the equivalence of (\ref{e22-3}) and (\ref{e22-4}), Xia et al. \cite{XIAZ2017} raised the following sufficient condition.
\begin{theorem}\label{th-2}{\rm\cite{XIAZ2017}}
Let $h_t(\x), t=1,\ldots,m_n$ denote polynomials on $\x\in\mathbb{R}^n$ of degree $d>1$, and $L=\{\x\in\mathbb{R}^n_+ : \a_i^\top\x-b_i=0, i=1,\ldots,m_l\}$. If the following conditions are satisfied, then (\ref{e22-3}) and (\ref{e22-4}) are equivalent:

\noindent(i) $h_t(\x)\geq 0$ for all $\x\in L$, where $t=1,\ldots,m_n$,

\noindent(ii) $\tilde{q}(\x)\geq 0$ for all $\x\in\{\x\in\mathbb{R}^n_+ : \tilde{l}_i(\x)=0, \tilde{h}_t(\x)=0, i=1,\ldots,m_l, t=1,\ldots,m_n\}$.
\end{theorem}

The primary motivation of this article is to establish an alternative characterization of the equivalence between POPs and CPTPs, without the need for additional verification. Compared to existing approaches in \cite{pena14,XIAZ2017}, our method requires only the recession cone of the feasible set of the original problem. This leads to a natural lifted convex formulation of the POP (see Section 3), through which an intrinsic connection between the solution of the original problem and that of the equivalent CPTP is revealed.

To apply the lifted convex reformulation and present an explicit equivalent CPTP, we focus on POPs with linear inequality constraints. It is worth noting that, for more general POPs, the CPTPs given in equations (\ref{e22-2}) and (\ref{e22-4}) are equivalent to the CPTP proposed in this paper, provided the sufficient conditions stated in Theorems \ref{th-1} and \ref{th-2} are satisfied. On one hand, if these conditions are not met, the equivalence between (\ref{e22-1}) and (\ref{e22-2}), or between (\ref{e22-3}) and (\ref{e22-4}), may fail to hold. Nevertheless, the CPTP model introduced in this work remains equivalent to the original problem (see Example \ref{exam-5}, or Section 4.3 of \cite{pena14}). On the other hand, a compelling and significant research direction lies in exploring the numerical trade-offs among various definitions such as those involving the recession cone. This issue is fundamentally shared by the authors of \cite{pena14,XIAZ2017}, and addressing it continues to be a meaningful avenue for future theoretical development.

We wish to clarify that the aim of this paper is not to reduce the computational complexity associated with the sufficient conditions in Theorems \ref{th-1} and \ref{th-2}. Rather, our primary objective is to construct an equivalent completely positive optimization formulation for specially structured polynomial optimization problems without invoking additional sufficient conditions. This approach provides potential benefits for advancing future research in numerical algorithm development.

\section{Lifted convex formulations of POP}\label{sec_unit}
In this section, we lift certain potentially non-convex polynomial optimization problems to equivalent convex formulations involving tensor variables. We then apply this lifting framework to a class of polynomial optimization problems subject to linear inequality constraints, which can be reformulated as equivalent completely positive optimization problems (see Section 4 for details).

\subsection{The lifted model for inhomogeneous POPs}
In this subsection, we examine the equivalent lifted convex formulation of inhomogeneous polynomial optimization problems. To proceed, we first introduce a useful mapping that characterizes the relationship between polynomials and symmetric tensors.

\begin{definition}\label{def33} Let $\mathbb{R}_d[\x]=\{p(\x)~|~{\rm deg}(p)\leq d\}$ denote the set of polynomials with dimension $n$ and degree at most $d$. For general polynomial $f(x)=\sum_{\alpha\in\mathbb{Z}^n_+} f_{\alpha}\x^{\alpha}\in \mathbb{R}_d[\x]$, the corresponding tensor
$\A_f\in S^{d,n+1}$ is defined as follows:
$$
(\A_f)_{\pi(i_1,i_2,\cdots, i_d)}=\frac{(d-|\alpha|)!\alpha_1!\cdots\alpha_n!}{d!}f_{\alpha},
$$
where $\pi(i_1,i_2,\cdots, i_d)$ denotes any permutation of $(i_1,i_2,\cdots, i_d)$, $\alpha$ is the unique exponent such that $\x^{\alpha}=x_1^{\alpha_1}x_2^{\alpha_2}\cdots x_n^{\alpha_n}=x_{i_1}x_{i_2}\cdots x_{i_d}$ and
$|\alpha|=\alpha_1+\alpha_2+\cdots+\alpha_n$.
\end{definition}

In Definition \ref{def33}, $\A_f$ is a symmetric tensor with entries defined by the coefficients of $f(\x)$. According to Definition \ref{def32}, we adopt $M_d(1,\x)=M_d((1,\x^\top)^\top)\in S^{d,n+1}$ for the sake of simplicity. Then,
it follows that
$$
f(\x)=\langle \A_f, M_d(1,\x)\rangle.
$$
For the benefit of readers,  we give an example to explain the detail of the equation above.
\begin{example} Suppose that $f(\x)\in\mathbb{R}_3[\x]$ is a polynomial function such that
$$
f(\x)=x_1^3+x_2^3+x_3^3+x_1x_2+x_2x_3+x_1+x_2-3.
$$
Let $\A=\A_f\in S^{3,4}$ and $\mathcal{M}=M_d(1,\x)\in S^{3,4}$. Then,  for any $\x\in\mathbb{R}^3$, it holds that $f(\x)=\langle \A, \mathcal{M}\rangle$, where
$$
\A_{0,:,:}=\left(
\begin{matrix}
-3 & \frac{1}{3} & \frac{1}{3} & 0\\
\frac{1}{3} & 0 & \frac{1}{6} & 0\\
\frac{1}{3} & \frac{1}{6} & 0 & \frac{1}{6}\\
0 & 0 &  \frac{1}{6} & 0
\end{matrix}
\right),
\A_{1,:,:}=\left(\begin{matrix}
 \frac{1}{3}& 0& \frac{1}{6} &0\\
0& 1& 0&0\\
 \frac{1}{6}& 0& 0&0\\
0&0 &0 &0
\end{matrix}
\right),
\A_{2,:,:}=\left(\begin{matrix}
\frac{1}{3}&\frac{1}{6} &0 &\frac{1}{6}\\
\frac{1}{6}&0 &0 &0\\
0&0 &1 &0\\
\frac{1}{6}& 0&0 &0
\end{matrix}
\right),
\A_{3,:,:}=\left(\begin{matrix}
0& 0&\frac{1}{6} &0\\
0&0 &0 &0\\
\frac{1}{6}& 0&0 &0\\
0&0 &0 &1
\end{matrix}
\right),
$$
and
$$
\begin{aligned}
&\mathcal{M}_{0,:,:}=\left(\begin{matrix}
1&x_1 & x_2&x_3\\
x_1&x_1^2 & x_1x_2&x_1x_3\\
x_2&x_2x_1 & x_2^2&x_2x_3\\
x_3& x_3x_1& x_3x_2&x_3^2
\end{matrix}
\right),~~~~~~~~~~
\mathcal{M}_{1,:,:}=\left(\begin{matrix}
x_1& x_1^2& x_1x_2 &x_1x_3\\
x^2_1&x_1^3 & x_1^2x_2&x_1^2x_3\\
x_1x_2& x_1^2x_2& x_1x_2^2&x_1x_2x_3\\
x_1x_3&x_1^2x_3 &x_1x_2x_3 &x_1x_3^2
\end{matrix}
\right),\\
&\mathcal{M}_{2,:,:}=\left(\begin{matrix}
x_2&x_2x_1 &x_2^2 &x_2x_3\\
x_2x_1&x_2x_1^2 &x_2^2x_1 &x_2x_1x_3\\
x_2^2&x_2^2x_1 &x_2^3 &x_2^2x_3\\
x_2x_3&x_2x_3x_1&x_2^2x_3 &x_2x_3^2
\end{matrix}
\right),~~
\mathcal{M}_{3,:,:}=\left(\begin{matrix}
x_3& x_3x_1&x_3x_2 &x_3^2\\
x_3x_1&x_3x_1^2 &x_3x_1x_2 &x_3^2x_1\\
x_3x_2& x_3x_2x_1&x_3x_2^2 &x_3^2x_2\\
x_3^2&x_3^2x_1 &x_3^2x_2 &x_3^3
\end{matrix}
\right).
\end{aligned}
$$
\end{example}\medskip

For inhomogeneous polynomial $f(\x)=\A_f\x^d\in\mathbb{R}_d[\x]$, it can be written as below:
$$
f(\x)=\A_d\x^d+\A_{d-1}\x^{d-1}+\cdots+\A_1\x+\A_0,
$$
where $\A_i\in S^{i,n}$. Particularly, $\A_0$ is a zero order tensor (scalar), $\A_1\in\mathbb{R}^n$ is a first order tensor (vector) and $\A_1\x$ means $\A_1^\top\x$. Let $\tilde{f}(\x)$ denote the homogeneous component of $f(\x)$ with the highest degree, which means that
$$
\tilde{f}(\x)=\A_d\x^d=\langle\A_d, M_d(\x) \rangle=\langle \A_f, M_d(0, \x)\rangle.
$$

We now consider the following inhomogeneous polynomial optimization problem:
\begin{equation}\label{e35}
\min f(\x)~~~{\rm s.t.}~~~\x\in F,
\end{equation}
where $f(x)\in \mathbb{R}_d[\x]$ and $F\subseteq\mathbb{R}^n$ is a nonempty feasible set. Similar to (\ref{e32}) for the homogeneous case, we have the following
lifted convex optimization problem:
\begin{equation}\label{e36}
\min~ \langle\A_f, \X\rangle~~~{\rm s.t.}~~~\X\in C_{F}+R_{F},
\end{equation}
where
$$
C_{F}={\rm conv}\{M_d(1,\x)~|~\x\in F \},~~R_{F}={\rm conv}\{M_d(0,\x)~|~\x\in REC_{F} \}.
$$
\begin{lemma}\label{lema32} Assume that (\ref{e35}) has an optimal solution $\bar{\x}$. Then, for any $\x\in REC_{F}$, it holds that
$\A_d\x^d\geq 0$.
\end{lemma}
{\bf Proof.} We prove the conclusion by contradiction. Assume that there is a $\bar{\y}\in REC_{F}$ satisfying $\A_d\bar{\y}^d<0$.
By the definition of recession cone, we know that for any $\lambda\geq 0$, it follows that $\bar{\x}+\lambda\bar{\y}\in F$ and
$f(\bar{\x}+\lambda\bar{\y})\geq f(\bar{\x})$. However, we have
$$
\begin{aligned}
&f(\bar{\x}+\lambda\bar{\y})=\A_d(\bar{\x}+\lambda\bar{\y})^d+\A_{d-1}(\bar{\x}+\lambda\bar{\y})^{d-1}+\cdots+\A_1(\bar{\x}+\lambda\bar{\y})+\A_0\\
&=C_d^d\lambda^d\A_d\bar{\y}^d+C_d^{d-1}\lambda^{d-1}\A_d\bar{\x}\bar{\y}^{d-1}+\cdots+C_d^1\lambda\A_d\bar{\x}^{d-1}\bar{\y}+\A_d\bar{\x}^d+\cdots+\A_0\\
&\rightarrow -\infty,~~~{\rm as}~~\lambda\rightarrow \infty,
\end{aligned}
$$
which contradicts the conditions, and the desired result holds.
\qed

When the optimal values of (\ref{e35}) and (\ref{e36}) can be attained, let $\mathbb{S}_1, \mathbb{S}_2$
be optimal solution sets of (\ref{e35}) and (\ref{e36}), respectively. By Lemma \ref{lema32}, the following results holds.

\begin{theorem}\label{them-3} Assume the optimal value of (\ref{e35}) is attained.
Suppose $\bar{\x}\in\mathbb{S}_1$. Then, the problems (\ref{e35}) and (\ref{e36}) are equivalent
in the sense that they share the same optimal value. Furthermore, it holds that $\mathbb{S}_2={\rm conv}\{M_d(1,\x) : \x\in\mathbb{S}_1\}$.
\end{theorem}
{\bf Proof.} Since $\x\in\mathbb{S}_1$, we have the optimal value such that $f(\bar{\x})=\langle\A_f, M_d(1,\bar{\x})\rangle$, and $M_d(1,\bar{\x})$ is feasible for problem (\ref{e36}). This means that
the optimal value of (\ref{e35}) is greater than the optimal value of (\ref{e36}).
Suppose $\bar{\X}\in \mathbb{S}_2$ and assume that
$$
\bar{\X}=\sum_{i=1}^k\lambda_iM_d(1, \x_i)+\sum_{j=1}^l\mu_jM_d(0, \y_j),
$$
where $\lambda_i\geq 0, i\in[k], \sum_{i=1}^k\lambda_i=1$ and $\mu_j\geq 0, j\in[l]$. Then, it holds that
\begin{equation}\label{e37}
\begin{aligned}
\langle\A_f, \bar{\X}\rangle&=\sum_{i=1}^k\lambda_i\langle\A_f,M_d(1, \x_i)\rangle+\sum_{j=1}^l\mu_j\langle\A_f,M_d(0, \y_j)\rangle\\
&=\sum_{i=1}^k\lambda_i\langle\A_f,M_d(1, \x_i)\rangle+\sum_{j=1}^l\mu_j\A_d\y_j^d.
\end{aligned}
\end{equation}
By Lemma \ref{lema32}, it follows that
\begin{equation}\label{e38}
\begin{aligned}
\langle\A_f, \bar{\X}\rangle &=\sum_{i=1}^k\lambda_i\langle\A_f,M_d(1,\x_i)\rangle+\sum_{j=1}^l\mu_j\A_d\y_j^d\geq\sum_{i=1}^k\lambda_i\langle\A_f,M_d(1, \x_i)\rangle\\
&=\sum_{i=1}^k\lambda_if(\x_i)\geq\sum_{i=1}^k\lambda_if(\bar{\x})=f(\bar{\x})\sum_{i=1}^k\lambda_i=f(\bar{\x}),
\end{aligned}
\end{equation}
which implies that (\ref{e35}) and (\ref{e36}) have the same optimal value.

From (\ref{e37}) and (\ref{e38}) again, since $\langle\A_f, \bar{\X}\rangle=f(\bar{\x})$, we may choose $\mu_j=0, j\in[l]$ and take
$\bar{\X}=\sum_{i=1}^k\lambda_iM_d(1, \x_i),$ which implies that $\bar{\X}$ is a convex combination of $M_d(1,\x_i), i\in[k]$ and $\x_i\in\mathbb{S}_1$. Therefore, the desired results hold.
\qed

\begin{remark} The polynomial optimization problem (\ref{e35}) is not a convex problem when the feasible set $F$ is not a convex set. So, it is challenging to solve (\ref{e35}) directly. By the equivalence of (\ref{e35}) and (\ref{e36}), in Section 4, we will show for some special cases, the lifted convex formulations (\ref{e36}) can be reformulated into completely positive optimization problems.
Furthermore, it is shown that, under
some conditions, the dual problems of the original completely positive optimization problem is strictly feasible, which means that strong duality holds and existing numerical methods for
both primal and dual problems may be applied.
\end{remark}

To verify the performance of Theorem \ref{them-3}, we consider the following example \cite{pena14}.
\begin{example}\label{exam-5} Consider the problem
\begin{equation}\label{e39}
\begin{aligned}
&\min~ f=4x-y-2x^2-2xy-y^2\\
&~{\rm s.t.}~~ x^2-xy=0,\\
&~~~~~~~y^2-y=0,\\
&~~~~~~~x\geq0,~y\geq 0.\\
\end{aligned}
\end{equation}
It is not difficult to see that the problem has only three feasible points $F=\{(0, 0),(0,1), (1, 1)\}$ and the optimal solutions of (\ref{e39}) are $(0,1), (1, 1)$ with optimal value equal to $-2$.
By Definition \ref{def33}, we have the corresponding second order tensor $\A_f$ such that
$$
\A_f=\left(
\begin{matrix}
0&2 &-\frac{1}{2} \\
2&-2 &-1 \\
-\frac{1}{2}&-1 &-1
\end{matrix}
\right).
$$
Since the recession cone of $F$ including the unique zero vector, we have the lifted model:
\begin{equation}\label{e311}
\min~ \langle\A_f, \X\rangle~~~{\rm s.t.}~~~\X\in {\rm conv}\{M_2(1,\x)~|~\x\in F \},
\end{equation}
and the set of optimal solutions
of (12) is
$$
\mbox{conv}\{M_2(1, 0, 1), M_2(1, 1, 1)\} =\left\{\left[
\begin{array}{lll}
1& \alpha & 1 \\
\alpha & \alpha & \alpha \\
1& \alpha & 1
\end{array}
\right]:0\leq\alpha\leq 1 \right\}.
$$
with optimal value $-2$.
\end{example}

\subsection{The lifted model for homogeneous POPs}

In this section, we consider the lifted model for homogeneous POPs since that many practical problems are related with homogeneous polynomial optimization problems, which are not convex optimization problems.

Let $f(\x)\in \mathbb{R}_d[\x], \x\in\mathbb{R}^n$ be a homogeneous polynomial. Consider the following polynomial optimization problem:
\begin{equation}\label{e31}
\min f(\x)~~{\rm s.t.}~~ \x\in F,
\end{equation}
where $f(\x)=\A^h_f\x^d$, and $F\subseteq\mathbb{R}^n$ is a nonempty feasible set. The entries of tensor $\A^h_f$ are constructed by the coefficients of the polynomial. By (\ref{e21}),
we have the following lifted optimization problem:
\begin{equation}\label{e32}
\min~\langle\A^h_f, \mathcal{X}\rangle~~{\rm s.t.}~~\X\in C^h_F+R^h_F.
\end{equation}

To prove the equivalence between (\ref{e31}) and (\ref{e32}), as in the inhomogeneous case, let $\A_i, i\in\{0,1,\cdots,d-1\}$ be zero tensors, then it reduces to homogeneous polynomial. Let $\mathbb{S}_3, \mathbb{S}_4$ be optimal solution sets of (\ref{e31}) and (\ref{e32}) respectively.
By a similar proof of Theorem \ref{them-3}, we have the following result.

\begin{theorem}\label{them-4} Assume the optimal value of (\ref{e31}) is attained.
Suppose $\bar{\x}\in\mathbb{S}_3$. Then, the problems (\ref{e31}) and (\ref{e32}) are equivalent
in the sense that they share the same optimal value. Furthermore, it holds that $\mathbb{S}_4={\rm conv}\{M_d(1,\x) : \x\in\mathbb{S}_3\}$.
\end{theorem}

\begin{remark} The primary distinction between the original problem (\ref{e31}) and the lifted model (\ref{e32}) lies in their structural properties: while the original problem may involve a non-convex objective function defined over a non-convex feasible set, the equivalent formulation (\ref{e32}) constitutes a convex optimization problem.
Therefore, the optimal solution of (\ref{e32}) includes global optimal solutions of (\ref{e31}).
For some special cases such as $M_d(\x)\in C^h_F+R^h_F$ and $\x\in\mathbb{R}^n_+$, (\ref{e32}) is a completely positive
tensor program. We will study this case in next section.

It is important to observe that the construction of the lifted problem (\ref{e32}) does not aim to remove the challenges associated with the original formulation (\ref{e31}); instead, these difficulties are embedded within the feasible set of (\ref{e32}), as determined by the definition of the set $F$.

\end{remark}

By Theorem \ref{them-4} and the definition of recession cones, we have the following corollary.
\begin{corollary}\label{corol31}
Assume that the feasible set $F$ of (\ref{e31}) is bounded and the optimal solution is attained. Then,  (\ref{e31}) shares the same optimal value with
the following convex optimization problem:
$$
\min~\langle\A^h_f, \mathcal{X}\rangle~~{\rm s.t.}~~\X\in {\rm conv}\{M_d(\x)~|~\x\in F\}.
$$
Furthermore, any optimal solution of the convex optimization problem can be denoted as a convex
combination of symmetric rank-1 tensors $M_d(\x)$, where $\x\in\mathbb{S}_3$.
\end{corollary}
{\bf Proof.} Since $F$ is bounded, its recession cone has only a zero vector,
which implies the desired result holds.
\qed

\begin{remark} 
The problems (\ref{e31}) and (\ref{e32}) are still equivalent when the constraint of (\ref{e32}) is replaced by
\begin{equation}
\min~\langle\A_f, \mathcal{X}\rangle~~{\rm s.t.}~~\X\in C^h_F.
\end{equation}
However, we will keep the constraints of (\ref{e32}) because the model (\ref{e32}) has advantages for the construction of completely positive reformulation in some
typical cases no matter the original polynomial is homogeneous or inhomogeneous.
\end{remark}

Notice that many practical problems are related with homogeneous polynomial optimization problems, which are not convex optimization problems.
By Theorem \ref{them-4}, some typical examples can be lifted to convex problems with tensor variables.
\begin{example}\label{exam31} Suppose $\A\in S^{d,n}$ is a symmetric tensor. Its minimal Z-eigenvalue equals the optimal value of the following optimization problem:
$$
\min \A\x^d~~~{\rm s.t.}~~\|\x\|=1,~\x\in\mathbb{R}^n,
$$
where $\|\cdot\|$ denotes the Euclidean norm of the vector. By a direct computation, the equivalent convex optimization problem of the minimal Z-eigenvalue problem is
$$\min~ \langle\A, \X\rangle~~~{\rm s.t.}~~~ \X\in {\rm conv}\{M_d(\x)~|~\|\x\|=1,~\x\in\mathbb{R}^n\}.$$
\end{example}

When the Euclidean norm is replaced by the norm $\|\cdot\|_m$, then it reduces to the minimal H-eigenvalue problem for a symmetric tensor.

\begin{example}\label{exam32}
Assume that $\A\in S^{d,n}$ is a symmetric tensor. To verify the copositiveness of $\A$, it is equivalent to solving the following optimization problem:
$$
\min\A\x^d~~~{\rm s.t.}~~~\x\in\mathbb{R}^n_+.
$$
By a direct computation, we know that the recession cone of $\mathbb{R}^n_+$ is itself. Thus, its equivalent convex optimization problem
is a completely positive tensor optimization problem:
$$\min~ \langle\A, \X\rangle~~~{\rm s.t.}~~~ \X\in {\rm conv}\{M_d(\x)~|~\x\in\mathbb{R}_+^n\}.$$
\end{example}

\section{Equivalent CPTPs for general POPs with linear inequality constraints}\label{sec03}
We are now ready to state our main results. By the lifted convex models (\ref{e36}) and (\ref{e32}), we construct the equivalent CPTP of the POP with explicit linear inequality constraints.
Then, the dual problem of the CPTP i.e., the copositive optimization problem is given.
In order to guarantee strong duality for the primal-dual pair under consideration,
we show that the copositive optimization problem is strictly feasible.

In the following sections, we first consider problems (\ref{e31}) and (\ref{e35}) with constraints defined by the intersection of a finite number of half-spaces such as
\begin{equation}\label{e41}
F=\{\x\in\mathbb{R}^n~|~B\x\leq \b,~\x\geq\0 \},
\end{equation}
where $B\in\mathbb{R}^{m,n}$ and $\b\in\mathbb{R}^m$. By (\ref{e21}) and Theorem \ref{them-4}, it is sufficient to provide an equivalent characterization of
$C_{F}+R_{F}$ over the cone of completely positive tensors.

To present the main results of this section, we first note that an optimal solution may not necessarily exist in general. To the best of our knowledge, additional sufficient conditions need to be imposed. For instance, the objective function is bounded from below, the constraint set is compact, or the symmetric tensor corresponding to the objective polynomial is a positive semi-definite tensor etc.
In this paper, for the sake of convenience, we do not conduct our study for any specific scenario, and we always have attainment as an assumption.

\subsection{Homogeneous POP with constraint set (\ref{e41}) }
Let $f(\x)\in \mathbb{R}_d[\x], \x\in\mathbb{R}^n$ be a homogeneous polynomial. Consider the following polynomial optimization problem:
\begin{equation}\label{e42}
\min f(\x)~~{\rm s.t.}~~ \x\in F,
\end{equation}
where $f(\x)=\A^h_f\x^d$.

To establish the equivalent CPTP formulation, we proceed in two steps. This approach is inspired by the proof for quadratic optimization problems in \cite{Ng2021}, where the author investigated an equivalent completely positive optimization framework for quadratic polynomial problems.
To move on, we first present several useful notations. By the definition of $F$, there always exists a vector $\alpha\in\mathbb{R}^n_+$ such that $\alpha^\top\x\leq 1$ for
$\x\in F$ (The trivial choice is $\alpha=\0$).

Let $\a=(\alpha^\top, 1)^\top\in\mathbb{R}^{n+1}_+$. Define the matrix $A$ with entries:
$$
A=\left(\begin{matrix}
\b\a^\top-(B,\0)\\
I_{n+1}\\
\end{matrix}
\right)\in\mathbb{R}^{(m+n+1)\times(n+1)},
$$
where $I_{n+1}$ is the identity matrix. Assume that $\A^h_f=(a_{i_1i_2\cdots i_d})$ in (\ref{e42}). Define $\bar{\A}=(\bar{a}_{i_1i_2\cdots i_d})\in S^{m, n+1}$ satisfying
$$
\bar{a}_{i_1i_2\cdots i_d}=
\left\{
\begin{aligned}
&a_{i_1i_2\cdots i_d},&i_1,i_2,\cdots,i_d\in[n],\\
&~~~~0, &{\rm otherwise}.~~~~~~~~~~~
\end{aligned}
\right.
$$
Let $\bar{f}(\y)=\bar{\A}\y^d$. Then, we have the following polynomial optimization problem:
\begin{equation}\label{e43}
\min \bar{f}(\y)~~~~~{\rm s.t.}~~\y\in \overline{F}=\{\y\in\mathbb{R}^{n+1}~|~A\y\geq0,~\a^\top\y=1\}.
\end{equation}
Furthermore, we obtain the following result.
\begin{lemma}\label{lema41} If problems (\ref{e42}) and (\ref{e43}) have optimal solutions, then they share
the same optimal value.
\end{lemma}
{\bf Proof.} Assume $\bar{\x}\in F$ and $\bar{\y}\in\overline{F}$ are optimal solutions of (\ref{e42}) and (\ref{e43}),  respectively. It follows that
$$B\bar{\x}\leq \b,~\bar{\x}\geq\0,~\a^\top\y=\alpha^\top\bar{\x}+1-\alpha^\top\bar{\x}=1,$$
where
$$\y=(\bar{\x}^\top, 1-\alpha^\top\bar{\x})^\top\in\mathbb{R}^{n+1}.$$
Then, we have that
$$
A\y=\left(
\begin{matrix}
\b-B\bar{\x}\\
\y
\end{matrix}
\right)\geq \0,
$$
which means that $\y\in\overline{F}$. Therefore, we obtain that $
\bar{f}(\y)=f(\bar{\x})\geq \bar{f}(\bar{\y}),
$
and the optimal value of (\ref{e42}) is greater than the optimal value of (\ref{e43}).

On the other hand, $\bar{\y}$ can be written as two blocks such that $\bar{\y}^\top=(\tilde{\y}^\top, y_{n+1})$, where $\tilde{\y}^\top=(y_1,y_2,\cdots, y_n)\in\mathbb{R}^n$. Since $\bar{\y}\in \overline{F}$ and $A\bar{\y}\geq \0$, we obtain that
$$
B\tilde{\y}\leq\b,~~~\tilde{\y}\geq\0,
$$
which implies that $\tilde{\y}\in F$ and
$$
\A\tilde{\y}^d=\bar{\A}\bar{\y}^d\geq\A\bar{\x}^d,
$$
which means that the optimal value of (\ref{e43}) is greater than the optimal value of (\ref{e42}) and the desired results hold.
\qed
\begin{remark}\label{r1}
From the proof of Lemma \ref{lema41}, we know that $\bar{\x}$ is an optimal solution of (\ref{e42}) if and only if $\bar{\y}=(\bar{\x}^\top, 1-\alpha^\top\bar{\x})^\top\in\mathbb{R}^{n+1}$ is an optimal solution of (\ref{e43}).
\end{remark}

Combining Remark \ref{r1} with Lemma \ref{lema41}, it is sufficient to construct the equivalent CPTP of (\ref{e43}) instead of (\ref{e42}).
\begin{lemma}\label{lema42}
For optimization problem (\ref{e43}), we have the following results.

\noindent{\rm(1)} The recession cone $REC_{\overline{F}}=\{\y\in\mathbb{R}^{n+1}~|~A\y\geq\0,~\a^\top\y=0\}$.

\noindent{\rm(2)} For feasible set $\overline{F}$, it holds that
$$
C^h_{\overline{F}}+R^h_{\overline{F}}=\left\{\X\in\mathcal{C}^*_{n+1,d}~|~\langle\X,M_d(\a)\rangle=1,  \X_{\times 1}A_{\times 2}A_{\times 3}\cdots _{\times d}A\in\mathcal{C}_{m+n+1,d}^*\right\}.
$$
\end{lemma}
{\bf Proof.}
\noindent{\rm(1)} This statement can be proved directly by the definition of the recession cone.

\noindent{\rm (2)} For any $\X\in C^h_{\overline{F}}+R^h_{\overline{F}}$, assume that
$$
\X=\sum_{i=1}^r\lambda_iM_d(\x_i)+\sum_{j=1}^l\mu_jM_d(\y_j),
$$
where $\lambda_i\geq0$, $\mu_j\geq0, i\in[r], j\in[l]$ and $\sum_{i=1}^r\lambda_i=1$. Then,  $\X\in\mathcal{C}_{n+1,d}^*$. Since $\x_i\in \overline{F}$, $\y_j\in REC_{\overline{F}}$, we have that $A\x_i\geq\0$, $A\y_j\geq\0$, $\a^\top\x_i=1$, $\a^\top\y_j=0$ and
$$
\langle\X,~M_d(\a)\rangle=\sum_{i=1}^r\lambda_i(\a^\top\x_i)^d+\sum_{j=1}^l\mu_j(\a^\top\y_j)^d=\sum_{i=1}^r\lambda_i=1,
$$
and
$$
\X_{\times 1}A_{\times 2}A\cdots _{\times d}A=\sum_{i=1}^r\lambda_iM_d(A\x_i)+\sum_{j=1}^l\mu_jM_d(A\y_j)\in\mathcal{C}_{m+n+1,d}^*,
$$
which implies that
$$
C^h_{\overline{F}}+R^h_{\overline{F}}\subseteq\left\{\X\in\mathcal{C}^*_{n+1}~|~\langle\X,M_d(\a)\rangle=1,  \X_{\times 1}A_{\times 2}A_{\times 3}\cdots _{\times d}A\in\mathcal{C}_{m+n+1,d}^*\right\}.
$$

To show the reverse inclusion, for $\X=(\X_{i_1i_2\ldots i_d})\in\mathcal{C}^*_{n+1,d}$ satisfying
$$
\langle\X,M_d(\a)\rangle=1, ~~ \X_{\times 1}A_{\times 2}A_{\times 3}\cdots _{\times d}A\in\mathcal{C}_{m+n+1,d}^*,
$$
let $\widehat{\X}\in\mathcal{C}^*_{m+n+1,d}$ with entries such that
$$
\widehat{\X}_{i_1i_2\ldots i_d}=\left\{
\begin{array}{lll}
&\X_{i_1i_2\ldots i_d}, &~~ i_1,i_2,\ldots, i_d\in[n+1],\\
&~~~~~0, &~~{\rm otherwise}.\\
\end{array}
\right.
$$
Denote $\widehat{\a}=(\a^\top,\0^\top)^\top\in\mathbb{R}^{m+n+1}_+$. Then it is clear that
\begin{equation}\label{e222}
\langle\widehat{\X},M_d(\widehat{\a})\rangle=\langle\X,M_d(\a)\rangle=1.
\end{equation}
On the other hand, there exists matrix $D$ that can be appended to $A$ such that $\widehat{A}=(A,D)\in\mathbb{R}^{(m+n+1)\times(m+n+1)}$ is a nonsingular matrix.
By the construction, it holds that
$$
\widehat{\X}_{\times 1}\widehat{A}_{\times 2}\widehat{A}_{\times 3}\cdots _{\times d}\widehat{A}=\X_{\times 1}A_{\times 2}A_{\times 3}\cdots _{\times d}A\in\mathcal{C}_{m+n+1,d}^*.
$$
Therefore, there are vectors $\widehat{\z}_1,\widehat{\z}_2,\cdots,\widehat{\z}_k\in\mathbb{R}_+^{m+n+1}$ such that
$$
\widehat{\X}_{\times 1}\widehat{A}_{\times 2}\widehat{A}_{\times 3}\cdots _{\times d}\widehat{A}=\sum_{i=1}^kM_d(\widehat{\z}_i).
$$
Let $\widehat{\y}_i=\widehat{A}^{-1}\widehat{\z}_i$ for $i=1,\ldots,k$, where $\widehat{A}^{-1}$ denotes the inverse of $\widehat{A}$. Then, we have that
$$
\begin{aligned}
\widehat{\X}&=\left(\widehat{\X}_{\times 1}\widehat{A}_{\times 2}\widehat{A}_{\times 3}\cdots _{\times d}\widehat{A}\right){_{\times 1}\widehat{A}^{-1}_{\times 2}\widehat{A}^{-1}_{\times 3}\cdots _{\times d}\widehat{A}^{-1}}\\
&=\sum_{i=1}^kM_d(\widehat{A}^{-1}\widehat{\z}_i)=\sum_{i=1}^kM_d(\widehat{\y}_i).
\end{aligned}
$$
Since $\widehat{\X}_{ii\ldots i}=0$ for all $i=n+2,n+3,\ldots,m+n+1$, it follows that $\widehat{\y}_j=(\y_j^\top,\0^\top)^\top$, where $\y_j$ denotes the first $n+1$ components of $\widehat{\y}_j, j\in[k]$.
By the structure of $\widehat{\X}$ and $\X$, we know that
$$
\X=\sum_{i=1}^kM_d(\y_i)~~{\rm and}~~A\y_i=\widehat{A}\widehat{\y}_i=\widehat{\z}_i\geq\0.
$$
Combining this with the notion of
$$
A=\left(\begin{matrix}
\b\a^\top-(B,\0)\\
I_{n+1}\\
\end{matrix}
\right),
$$
it follows that $\y_i\geq\0, i\in[k]$.
By the fact that $\y_i\geq\0, \a\geq \0$, it means that $\y_i^\top\a\geq 0, i\in[k]$.
To continue, denote
\begin{equation}\label{e45}
K^+=\left\{i\in[k]~|~\a^\top\y_i>0\right\},~~~K^==\{i\in[k]~|~\a^\top\y_i=0\}.
\end{equation}
For any $i\in K^+$, let $\lambda_i=\a^\top\y_i$ and $\z_i=\frac{1}{\lambda_i}\y_i$.
By this notation, we obtain that
$$
\X=\sum_{i\in K^+}\lambda_i^dM_d(\z_i)+\sum_{j\in K^=}M_d(\y_j).
$$
Next, we will prove that $\X\in C^h_{\overline{F}}+R^h_{\overline{F}}$. By (\ref{e45}),
we know that $A\z_i=\frac{1}{\lambda_i}A\y_i\geq \0, i\in K^+$. Furthermore, by (\ref{e222}), it implies that
$$
\a^\top\z_i=\frac{1}{\lambda_i}\a^\top\y_i=1,~~\forall~i\in K^+, \sum_{i\in K^+}\lambda_i^m=1.
$$
Combining this with (\ref{e45}), we obtain that
$$
\sum_{i\in K^+}\lambda_i^dM_d(\z_i)\in C^h_{\overline{F}},~~\sum_{j\in K^=}M_d(\y_j)\in R^h_{\overline{F}},
$$
which means that
$$
C^h_{\overline{F}}+R^h_{\overline{F}}\supseteq\left\{\X\in\mathcal{C}^*_{n+1}~|~\langle\X,M_d(\a)\rangle=1,  \X_{\times 1}A_{\times 2}A_{\times 3}\cdots _{\times d}A\in\mathcal{C}_{m+n+1,d}^*\right\},
$$
and the desired results hold.
\qed

By Theorem \ref{them-4} and Lemmas \ref{lema41}-\ref{lema42}, we have the following results.
\begin{theorem}\label{them41} Assume (\ref{e42}) has an optimal solution. Then,  (\ref{e42}) is equivalent with the following completely positive optimization problem:
$$
\begin{aligned}
&\min~ \langle\A^h_f, \X\rangle\\
&~{\rm s.t.}~ \X\in\left\{\X\in\mathcal{C}^*_{n+1}~|~\langle\X,M_d(\a)\rangle=1,  \X_{\times 1}A_{\times 2}A_{\times 3}\cdots _{\times d}A\in\mathcal{C}_{m+n+1,d}^*\right\}.
\end{aligned}
$$
\end{theorem}

\subsection{Inhomogeneous POP with constraint set (\ref{e41})}
In this section, we consider the inhomogeneous polynomial optimization problem.
If $f(\x)\in\mathbb{R}_d[\x]$ and $\x\in\mathbb{R}^n$, assume that
$$
f(\x)=\A_d\x^d+\A_{d-1}\x^{d-1}+\cdots+\A_1\x+\A_0,
$$
where $\A_j=(a^j_{i_1i_2\cdots i_j})\in S^{j,n}$ and $\A_d$ is a tensor with at least one nonzero element. Consider the POP with constraint (\ref{e41}):
\begin{equation}\label{e46}
\min f(\x)~~~{\rm s.t.}~~\x\in F.
\end{equation}
To construct the equivalent completely optimization problem of (\ref{e46}), define the
symmetric tensor $\bar{\A}_j=(\bar{a}^j_{i_1i_2\cdots i_j})\in S^{j, n+1}$ such that
$$
\bar{a}^j_{i_1i_2\cdots i_j}=
\left\{
\begin{aligned}
&a^j_{i_1i_2\cdots i_j},~i_1,i_2,\cdots,i_j\in[n],\\
&0, ~~~~~~~~~~{\rm Otherwise}.
\end{aligned}
\right.
$$
By this, we have the polynomial $\bar{f}(\y), \y\in\mathbb{R}^{n+1}$ such that
$$
\bar{f}(\y)=\bar{\A}_d\y^d+\bar{\A}_{d-1}\y^{d-1}+\cdots+\bar{\A}_1\y+\bar{\A}_0.
$$
Then, there is another polynomial optimization problem such that
\begin{equation}\label{e47}
\min \bar{f}(\y)~~~~~{\rm s.t.}~~\y\in \overline{F}=\{\y\in\mathbb{R}^{n+1}~|~A\y\geq0,~\a^\top\y=1\},
\end{equation}
where the matrix $A$ and vector $\a$ are defined as in Subsection 4.1.
If $\y^\top=(\x^\top, y_{n+1})\in\mathbb{R}^{n+1}, \x\in\mathbb{R}^n$, it is obvious that $\bar{f}(\y)=f(\x)$. By a similar proof with Lemma \ref{lema41}, we have the following result.
\begin{lemma}\label{lema43} If problems (\ref{e46}) and (\ref{e47}) have optimal solutions, they share the same optimal value.
\end{lemma}

By Theorem \ref{them-3}, we know that the problem (\ref{e47}) is equivalent with the following convex problem:
\begin{equation}\label{e48}
\min~ \langle\A_{\bar{f}}, \X\rangle~~~{\rm s.t}~~\X\in C_{\overline{F}}+R_{\overline{F}},
\end{equation}
where
$$
C_{\overline{F}}={\rm conv}\{M_d(1, \y)~|~\y\in\overline{F}\},~~R_{\overline{F}}={\rm conv}\{M_d(0, \y)~|~\y\in REC_{\overline{F}}\}.
$$

In order to construct the equivalent CPTP of (\ref{e46}), define
$\bar{A}\in\mathbb{R}^{t\times(n+2)}$, and let $\bar{\A}\in S^{d,n+2}$ be a symmetric rank-1 tensor such that
$$
\bar{A}=
\left(
\begin{aligned}
&1, ~~-\a^\top\\
&1, ~~-\a^\top\\
&\vdots~~~~~~~\vdots \\
&1, ~~-\a^\top
\end{aligned}
\right),~~~~
\bar{\A}=
\left(
\begin{matrix}
1\\
-\a
\end{matrix}
\right)\circ\left(
\begin{matrix}
1\\
-\a
\end{matrix}
\right)\circ\cdots\circ\left(
\begin{matrix}
1\\
-\a
\end{matrix}
\right),
$$
where $t$ is a positive integer. Moreover, it is necessary to present another simple abbreviations. For $\x=(x_1,x_2,\cdots,x_{n+2})^\top\in\mathbb{R}^{n+2}$ and $\X=M_d(\x)\in S^{d,n+2}$, denote $\widetilde{\x}\in\mathbb{R}^{n+1}$ and $\widetilde{\X}\in S^{d, n+1}$ such that
$$
\widetilde{\x}=(x_2,\cdots,x_{n+2})^\top,~~\widetilde{\X}=\widetilde{\x}\circ\widetilde{\x}\circ\cdots\circ\widetilde{\x}.
$$
\begin{lemma}\label{lema44}
Assume that $C_{\overline{F}}, R_{\overline{F}}$ are defined as in (\ref{e48}). Then, $C_{\overline{F}}+R_{\overline{F}}$ is equivalent with the following set:

\noindent$
\left\{\X\in\mathcal{C}^*_{n+2}:\langle\widetilde{\X},~ M_d(\a)\rangle=1, ~\langle\X,\bar{\A}\rangle=0, ~\widetilde{\X}_{\times 1}A_{\times 2}\cdots_{\times d}A\in\mathcal{C}_{m+n+1,d}^*,\right.$

$\left.\X_{\times 1}\bar{A}_{\times 2}\cdots_{\times d}\bar{A}\in\mathcal{C}_{t,d}^*\right\}.
$
\end{lemma}
{\bf Proof.}
For any $\X=(\X_{i_1i_2\ldots i_d})\in\mathcal{C}^*_{n+2}$ satisfying that
\begin{equation}\label{e49}
\left\{
\begin{aligned}
&\langle\widetilde{\X}, M_d(\a)\rangle=1,~~ \langle\X,\bar{\A}\rangle=0, \\
&\widetilde{\X}_{\times 1}A_{\times 2}A\cdots_{\times d}A\in\mathcal{C}_{m+n+1,d}^*,~~\X_{\times 1}\bar{A}_{\times 2}\bar{A}\cdots_{\times d}\bar{A}\in\mathcal{C}_{t,d}^*,
\end{aligned}
\right.
\end{equation}
define $\widehat{\X}\in\mathcal{C}^*_{m+n+2,d}$ with entries such that
$$
\widehat{\X}_{i_1i_2\ldots i_d}=\left\{
\begin{array}{lll}
&\X_{i_1i_2\ldots i_d}, &~~ i_1,i_2,\ldots, i_d\in[n+2],\\
&~~~~~0, &~~{\rm otherwise}.\\
\end{array}
\right.
$$
Similarly, define tensor $\widehat{\A}\in\mathcal{\C}^*_{m+n+2,d}$ and vector $\widehat{\a}\in\mathbb{R}^{m+n+1}$ as below
$$
\widehat{\A}=
\left(
\begin{matrix}
1\\
-\a\\
\0
\end{matrix}
\right)\circ\left(
\begin{matrix}
1\\
-\a\\
\0
\end{matrix}
\right)\circ\cdots\circ\left(
\begin{matrix}
1\\
-\a\\
\0
\end{matrix}
\right),~~\widehat{\a}=(\a^\top,\0^\top)^\top\in\mathbb{R}^{m+n+1}.
$$
Clearly, by the constructions, it holds that
\begin{equation}\label{e333}
\langle\widetilde{\widehat{\X}}, M_d(\widehat{\a})\rangle=\langle\widetilde{\X}, M_d(\a)\rangle=1~~{\rm and}~~\langle\widehat{\X},\widehat{\A}\rangle=\langle\X,\bar{\A}\rangle=0.
\end{equation}
Let matrix $\widehat{A}$ be defined as in Lemma \ref{lema42}, and $\widehat{\bar{A}}=(\bar{A},\0)\in\mathbb{R}^{t\times(m+n+2)}$. By a direct computation, we know that
$$
\widetilde{\widehat{\X}}_{\times 1}\widehat{A}_{\times 2}\cdots_{\times d}\widehat{A}=\widetilde{\X}_{\times 1}A_{\times 2}\cdots_{\times d}A\in\mathcal{C}_{m+n+1,d}^*
$$
and
\begin{equation}\label{e333-1}
\widehat{\X}_{\times 1}\widehat{\bar{A}}_{\times 2}\widehat{\bar{A}}_{\times 3}\cdots_{\times d}\widehat{\bar{A}}=\X_{\times 1}\bar{A}_{\times 2}\bar{A}_{\times 3}\cdots_{\times d}\bar{A}\in\mathcal{C}_{t,d}^*.
\end{equation}
Since $\widetilde{\widehat{\X}}_{\times 1}\widehat{A}_{\times 2}\cdots_{\times d}\widehat{A}\in\mathcal{C}_{m+n+1,d}^*$, there are vectors $\widehat{\z}_1,\widehat{\z}_2,\cdots,\widehat{\z}_k\in\mathbb{R}_+^{m+n+2}$ such that
$$
\widetilde{\widehat{\X}}_{\times 1}\widehat{A}_{\times 2}\widehat{A}_{\times 3}\cdots _{\times d}\widehat{A}=\sum_{i=1}^kM_d(\widetilde{\widehat{\z}}_i).
$$
Define $\widehat{\y}_i\in\mathbb{R}^{m+n+2}$ and $\widetilde{\widehat{\y}}_i=\widehat{A}^{-1}\widetilde{\widehat{\z}}_i$ for $i=1,2,\cdots,k$.
Then, it follows that $\widehat{A}\widetilde{\widehat{\y}}_i=\widetilde{\widehat{\z}}_i\geq\0$ and
\begin{equation}\label{e444}
\begin{aligned}
\widetilde{\widehat{\X}}&=\left(\widetilde{\widehat{\X}}_{\times 1}\widehat{A}_{\times 2}\cdots_{\times d}\widehat{A}\right){_{\times 1}\widehat{A}^{-1}_{\times 2}\cdots_{\times d}\widehat{A}^{-1}}\\
&=\left(\sum_{i=1}^kM_d(\widetilde{\widehat{\z}}_i)\right){_{\times 1}\widehat{A}^{-1}_{\times 2}\cdots_{\times d}\widehat{A}^{-1}}=\sum_{i=1}^kM_d(\widehat{A}^{-1}\widetilde{\widehat{\z}}_i)=\sum_{i=1}^kM_d\widetilde{\widehat{\y}}_i.
\end{aligned}
\end{equation}
Since $\widehat{\X}_{ii\cdots i}=0$ for all $i=n+3,n+4,\cdots,m+n+2$, it follows that $\widehat{\y}_j=(\y_j^\top,\0^\top)^\top$, where $\y_j$ denotes the first $n+2$ components of $\widehat{\y}_j, j\in[k]$.
Combining this with (\ref{e444}), we have $\X=\sum_{i=1}^kM_d(\y_i)$.

On the other hand, by (\ref{e333})-(\ref{e333-1}), it implies that $(\y_i)_1=\a^\top\widetilde{\y}_i$ for all $i\in[k]$.
From the fact that $\a\geq\0, \y_i\geq\0$, we have $\a^\top\widetilde{\y}_i\geq 0$. Denote
\begin{equation}\label{e410}
L^+=\{i\in[k]~|~\a^\top\widetilde{\y}_i>0\},~~~L^==\{j\in[k]~|~\a^\top\widetilde{\y}_j=0\}.
\end{equation}
For $i\in L^+$, let $\lambda_i=\a^\top\widetilde{\y}_i$ and $\z_i=\frac{1}{\lambda_i}\y_i$.
Then $\X$ can be denoted equivalently such that
$$
\X=\sum_{i\in L^+}\lambda_i^dM_d(\z_i)+\sum_{j\in L^=}M_d(\y_i).
$$
Next, we will prove that $\sum_{i\in L^+}\lambda_i^dM_d(\z_i)\in C_{\overline{F}}$ and $\sum_{j\in L^=}M_d(\y_i)\in R_{\overline{F}}$.
By (\ref{e410}), we know that $\a^\top\widetilde{\z}_i=1$, $(\z_i)_1=\frac{1}{\lambda_i}(\y_i)_1=\frac{1}{\lambda_i}\a^\top\widetilde{\y}_i=1, i\in L^+$ and $(\y_j)_1=\a^\top\widetilde{\y}_j=0, j\in L^=$.
Furthermore, by (\ref{e333}) again, we obtain that
\begin{equation}\label{e411}
\begin{aligned}
\langle\widetilde{\X},~ M_d(\a)\rangle&=\sum_{i\in L^+}\lambda_i^d(\a^\top\widetilde{\z}_i)^d+\sum_{j\in L^=}(\a^\top
\widetilde{\y}_j)^d=\sum_{i\in L^+}\lambda_i^d(\a^\top\widetilde{\z}_i)^d=\sum_{i\in L^+}\lambda_i^d=1.
\end{aligned}
\end{equation}
Combining (\ref{e411}) with (\ref{e49}) and (\ref{e410}), it follows that
$\X\in C_{\overline{F}}+R_{\overline{F}}$.

To prove the reverse inclusion, for $\X\in C_{\overline{F}}+R_{\overline{F}}$, suppose that
$$
\X=\sum_{i=1}^r\lambda_iM_d(1,\y_i)+\sum_{j=1}^l\mu_jM_d(0,\z_j),~\y_i\in\overline{F},~ \z_j\in RES_{\overline{F}},
$$
where $\lambda_i\geq 0, \mu_j\geq 0$ and $\sum_{i=1}^r\lambda_i=1$.
From definitions of $\overline{F}$ and $RES_{\overline{F}}$,
it's obvious that
$\X\in\mathcal{C}_{n+2,d}^*$, $\a^\top\y_i=1, i\in[r]$, $\a^\top\z_j=0, j\in[l]$ and, by a direct computation, it holds that
$$
\langle\X,~\tilde{\A}\rangle=0,~ \widetilde{\X}_{\times 1}A_{\times 2}A_{\times 3}\cdots_{\times d}A\in\mathcal{C}_{m+n+1}^*,~\X_{\times 1}\bar{A}_{\times 2}\bar{A}_{\times 3}\cdots_{\times d}\bar{A}\in\mathcal{C}_{t,d}^*.
$$
Moreover, we have that
$$
\langle\widetilde{\X}, ~M_d(\a)\rangle=\sum_{i=1}^r\lambda_i(\a^\top\y_i)^d+\sum_{j=1}^l\mu_j(\a^\top\z_j)^d=\sum_{i=1}^r\lambda_i=1.
$$
Therefore, the desired results hold.
\qed

From Theorem \ref{them-3} and Lemmas \ref{lema43}-\ref{lema44}, we have the following theorem.
\begin{theorem}\label{them42} For general polynomial optimization problem (\ref{e46}),
if it has an optimal solution. Then,  (\ref{e46}) is equivalent with the following completely positive program:
$$
\begin{aligned}
&\min~ \langle\A_{\bar{f}},~ \X\rangle\\
&~{\rm s.t.}~~ \langle\widetilde{\X}, ~M_d(\a)\rangle=1,~ \langle\X,~\bar{\A}\rangle=0,~ \X_{\times 1}\bar{A}_{\times 2}\bar{A}\cdots_{\times d}\bar{A}\in\mathcal{C}_{t,d}^*,\\
&~~~~~~~\widetilde{\X}_{\times 1}A_{\times 2}A\cdots_{\times d}A\in\mathcal{C}_{m+n+1,d}^*,~~ \X\in\mathcal{C}^*_{n+2,d}.
\end{aligned}
$$
\end{theorem}

\begin{remark}\label{remark41} ~~\par
\noindent({\rm 1}) In the definition of $\bar{A}$, the number $t$ can be any positive integer such that $t\geq1$.
Here, the role of $\bar{A}$ in $\X_{\times 1}\bar{A}_{\times 2}\bar{A}_{\times 3}\cdots_{\times d}\bar{A}\in\mathcal{C}_{t,d}^*$
is important and it makes all generating vectors $\y_i$ of $\X=\sum_{i=1}^kM_d(\y_i)$ satisfying $(\y_i)_1-\a^\top\widetilde{\y}_i\geq0$.

\noindent({\rm 2}) For the optimization problem in Theorem \ref{them42}, we can uniformly consider the variable tensor $\X$. For any $\X\in\mathcal{C}^*_{n+2,d}$. It holds that $\langle\widetilde{\X}, ~M_d(\a)\rangle=1$ is equivalent with
$\langle\X,~M_d(\bar{\a})\rangle=1$, where $\bar{\a}=(0, \a^\top)^\top\in\mathbb{R}^{n+2}$.
Similarly, $\widetilde{\X}_{\times 1}A_{\times 2}A_{\times 3}\cdots_{\times d}A\in\mathcal{C}_{m+n+1,d}^*$ is equivalent with
$\X_{\times 1}A'_{\times 2}A'_{\times 3}\cdots_{\times d}A'\in\mathcal{C}_{m+n+1,d}^*$, where  $A'=(\0, A)\in\mathbb{R}^{(m+n+1)\times(n+2)}$.
\end{remark}

\subsection{Equivalent copositive optimization problems}
By the duality of the completely positive tensor cone and the copositive tensor cone, we consider the dual optimization problems of (\ref{e42}) and (\ref{e46}) i.e., the copositive optimization problems.
To the best of our knowledge, the copositive optimization problem has many advantages in designing algorithms for quadratic optimization problems and high order POPs \cite{Bomze2000,Burer2009,MDFR21,WCYZ2024}.
Therefore, to guarantee the strong duality for the primal-dual pair under consideration, we show that the copositive optimization problem is strictly feasible.
Recall the equivalent model in Theorem \ref{them41}:
\begin{equation}\label{e412}
\begin{aligned}
&\min~ \langle\A^h_f,~ \X\rangle\\
&~{\rm s.t.}~ \X\in\left\{\X\in\mathcal{C}^*_{n+1}~|~\langle\X,~M_d(\a)\rangle=1,  \X_{\times 1}A_{\times 2}A_{\times 3}\cdots _{\times d}A\in\mathcal{C}_{m+n+1,d}^*\right\}.
\end{aligned}
\end{equation}
Its Lagrangian function is
$$
L(\X,\lambda,\mathcal{U},\mathcal{V})=\langle\A^h_f,~ \X\rangle-\lambda(\langle\X,~ M_d(\a)\rangle-1)+\langle\mathcal{U},~ \X_{\times 1}A_{\times 2}A_{\times 3}\cdots_{\times d}A\rangle-\langle\mathcal{V},~ \X\rangle,
$$
where $\mathcal{U}\in \mathcal{C}_{m+n+1,d}, \mathcal{V}\in\mathcal{C}_{n+1,d}$ are copositive tensors and $\lambda\in\mathbb{R}$.
It can be rewritten equivalently as
$$
L(\X,\lambda,\mathcal{U},\mathcal{V})=\left\langle\A^h_f-\lambda M_d(\a)+\mathcal{U}_{\times 1}A^\top_{\times 2}A^\top_{\times 3}\cdots_{\times d}A^\top-\mathcal{V},~ \X\right\rangle+\lambda.
$$
Therefore, the dual problem of (\ref{e412}) is the following copositive optimization problem:
\begin{equation}\label{e413}
\begin{aligned}
&\max~ \lambda\\
&~{\rm s.t.}~~ \A^h_f-\lambda M_d(\a)+\mathcal{U}_{\times 1}A^\top_{\times 2}A^\top_{\times 3}\cdots_{\times d}A^\top\in\mathcal{C}_{n+1,d},\\
&~~~~~~~\mathcal{U}\in\mathcal{C}_{m+n+1,d},~\lambda\in\mathbb{R}.
\end{aligned}
\end{equation}
Obviously, when $\mathcal{U}$ is a zero tensor and $\lambda$ is small enough, there exists $\lambda<0$ such that $\A_f-\lambda M_d(\a)$
is strictly copositive, which means that (\ref{e413}) is strictly feasible and the strong duality holds for (\ref{e412}) and (\ref{e413}).

Similarly, by Theorem \ref{them42} and Remark \ref{remark41}, the equivalent completely positive optimization of (\ref{e46}) can be rewritten as follows:
\begin{equation}\label{e414}
\begin{aligned}
&\min~ \langle\A_{\bar{f}},~ \X\rangle\\
&~{\rm s.t.}~~ \langle\X, ~M_d(\bar{\a})\rangle=1,~ \langle\X,~\tilde{\A}\rangle=0,~ \X_{\times 1}\bar{A}_{\times 2}\bar{A}_{\times 3}\cdots_{\times d}\bar{A}\in\mathcal{C}_{t,d}^*,\\
&~~~~~~~\X_{\times 1}A'_{\times 2}A'_{\times 3}\cdots_{\times d}A'\in\mathcal{C}_{m+n+1}^*,~~ \X\in\mathcal{C}^*_{n+2,d}.
\end{aligned}
\end{equation}
By a direct computation, we obtain the dual problem of (\ref{e414}) is
\begin{equation}\label{e415}
\begin{aligned}
&\max~ \lambda\\
&~{\rm s.t.}~~ \A_{\bar{f}}-\lambda M_d(\bar{\a})-\mu\tilde{\A}+\mathcal{U}^1_{\times 1}\bar{A}^\top_{\times 2}\cdots_{\times d}\bar{A}^\top+\mathcal{U}^2_{\times 1}A'^\top_{\times 2}\cdots_{\times d}A'^\top\in\mathcal{C}_{n+2,d},\\
&~~~~~~~\mathcal{U}^1\in\mathcal{C}_{t},~\mathcal{U}^2\in\mathcal{C}_{m+n+1,d},~\lambda,~\mu\in\mathbb{R}.
\end{aligned}
\end{equation}
It is obvious that the (\ref{e415}) is strictly feasible and the strong duality holds for (\ref{e414}) and (\ref{e415}).

\section{Conclusions}\label{sec_con}
In this paper, we investigate the equivalent reformulations of polynomial optimization problems (POPs) as completely positive tensor programs (CPTPs). Rather than assuming nonnegativity of the polynomial constraints over $\mathbb{R}^n_+$,  we consider a relaxed condition, imposing nonnegativity solely over the feasible region defined by linear inequality constraints. To facilitate this, we first introduce a general lifting framework that transforms POPs into convex programs involving tensor variables. This lifting approach is then employed to provide a broad characterization of POPs with linear constraints that can be equivalently represented as CPTPs. In addition, the duals of the reformulated CPTPs, namely copositive programs, are constructed and shown to be strictly feasible, thus ensuring that the strong duality holds.

To the best of our knowledge, unlike quadratic optimization problems, efficient numerical algorithms for completely positive tensor optimization problems remain unavailable due to the inherent complexity of these models. As is well known, the dual of a completely positive tensor optimization problem corresponds to a copositive tensor optimization problem, a relatively new area in mathematical optimization that generalizes copositive matrix optimization. While significant progress has been made in developing numerical algorithms for copositive matrix optimization (see \cite{MDFR21}), comparable advancements for copositive tensor optimization are still lacking. Investigating numerical methods for high-order copositive tensor optimization problems presents an interesting and promising direction for future research.

\medskip
\begin{acknowledgements}
The authors gratefully acknowledge the editor and two anonymous referees for their constructive comments, which significantly contributed to improving the quality of this paper.
This work was supported by Natural Science Foundation of China (12071249), Shandong Provincial Natural Science Foundation (ZR2024MA003,ZR2021JQ01), Hong Kong Innovation and Technology Commission (InnoHK Project CIMDA) and Hong Kong Research Grants Council (Project CityU 11204821).
\end{acknowledgements}

\noindent{\bf Availability of data and materials:} The data used to support the findings of this study are
available from the corresponding author upon request.

\end{document}